\documentclass{amsart}

\usepackage[dvipsnames]{xcolor}
\usepackage{tikz}
\usepackage{amsmath,bm,bbm,amsthm, amssymb}
\usepackage{mathtools}
\usepackage{fullpage}
\usepackage{color}
\usepackage{dsfont}
\usepackage{animate}
\usepackage{etoolbox}
\usepackage{comment}
\usepackage{pgf}
\usetikzlibrary{calc}
\usetikzlibrary{patterns}
\usetikzlibrary{arrows}
\usetikzlibrary{decorations.pathreplacing}
\usepackage[utf8]{inputenc}
\usepackage{pgfplots}

\usepackage{tcolorbox}
\usepackage{adjustbox}
\usepackage[final]{pdfpages}
\usepackage[ruled,vlined,linesnumbered,algosection,resetcount]{algorithm2e}
\usepackage{blkarray}
\usepackage{arydshln}
\usepackage{wrapfig}

\usepackage[textwidth=2cm]{todonotes}
\allowdisplaybreaks
\theoremstyle{plain}
\newtheorem{theorem}{Theorem}
\newtheorem{proposition}[theorem]{Proposition}
\newtheorem{corollary}[theorem]{Corollary}
\newtheorem{lemma}[theorem]{Lemma}

\theoremstyle{definition}

\theoremstyle{remark}
\newtheorem{remark}[theorem]{Remark}

\def\R{\mathbb R}

\def\P{\mathbb P}

\def\E{\mathbb E}

\def\eps{\varepsilon}

\def\xx{\mathbf x}
\def\yy{\mathbf y}

\def\00{\mathbf 0}

\def\Ga{\Gamma}

\def\cum{\ms{cum}}
\def\bet{\begin{theorem}}
\def\ent{\end{theorem}}
\def\bec{\begin{corollary}}
\def\enc{\end{corollary}}
\def\bep{\begin{proof}}
\def\enp{\end{proof}}
\def\f{\frac}

\def\a{\tau}
\def\g{\gamma}
\def\la{\lambda}

\def\su{\subseteq}
\def\ms{\mathsf}
\def\co{\colon}
\def\N{\mathbb N}
\def\mc{ \mathcal}
\def\ff{\infty}
\def\PP{\mc P}

\def\one{\mathds1}

\def\d{\mathrm d}

\def\r{\rho}

\renewcommand\le{\leqslant}
\renewcommand\ge{\geqslant}

\def\e{\varepsilon}
\def\bel{\begin{lemma}}
\def\enl{\end{lemma}}
\def\im{\item}
\def\been{\begin{enumerate}}
\def\enen{\end{enumerate}}
\def\pa{\partial}
\def\sm{\setminus}

\def\k_d{\kappa_d}

\def\s{\sigma}
\def\co{\colon}
\def\ti{\times}
\def\k{\kappa}

\def\bepr{\begin{proposition}}
\def\enpr{\end{proposition}}
\def\vp{\varphi}

\def\b{\beta}
\def\Var{\ms{Var}}

\def\De{\Delta}

\def\cuo{c_{U, 1}}
\def\cut{c_{U, 2}}
\def\cpo{c_{\vp, 1}}
\def\cpt{c_{\vp, 2}}

\begin{document}

\title{Cumulant method for weighted random connection models}

\author{Nils Heerten}
\author{Christian Hirsch}
\author{Moritz Otto}
\address[Nils Heerten]{Ruhr University Bochum, Bochum, Germany}
\email{nils.heerten@rub.de}
\address[Christian Hirsch, Moritz Otto]{Department of Mathematics, Aarhus University, Ny Munkegade 118, 8000 Aarhus C, Denmark}
\email{hirsch@math.au.dk, otto@math.au.dk}
\address[Christian Hirsch]{DIGIT Center, Aarhus University, Finlandsgade 22, 8200 Aarhus N, Denmark}


\begin{abstract}
	In this paper, we derive cumulant bounds for subgraph counts and power-weighted edge length in a class of spatial random networks known as weighted random connection models. This involves dealing with long-range spatial correlations induced by the profile function and the weight distribution. We start by deriving the bounds for the classical case of a Poisson vertex set, and then provide extensions to $\alpha$-determinantal processes. 
\end{abstract}
\maketitle
\vspace{0.3cm}


%
%
\section{Introduction}
\label{sec:int}

%
%
In the random geometric graph, the vertex set is a random set of points in Euclidean space, and any two vertices are connected by an edge if their distance is smaller than a fixed threshold. Such a random geometric graph is a model in a wide variety of applications such as wireless networks. Therefore, its asymptotic behavior for a large number of nodes has been studied extensively in the literature \cite{penrose}.  However, with the rise of network science, it has become apparent that often this model is not flexible enough to capture more complex structures. Driven by this need, the class of \emph{weighted random connection models (WRCM)} has recently appeared as a promising approach for describing phenomena in large complex networks \cite{glm2,komjathy2,komjathy}. Instead of a fixed connection threshold, the WRCM quantify the connectivity via a \emph{profile function} and a \emph{weight distribution}.  However, due to the complexity of the models, the analysis of their asymptotic behavior is still in its infancy. In this paper, we study the question of asymptotic normality of suitable functionals on the WRCM. This is relevant in order to put the analysis of data via WRCM  on a firm statistical foundation. 

%
%
One of the main tools in the analysis of asymptotic normality is the method of cumulants \cite{doring2021,saulis1991}. Its popularity stems from the level of generality where it is applicable and the variety of implications it draws. Moreover, one obtains specific non-asymptotic bounds for the rate of convergence to the normal distribution. However, the challenge in applying the method of cumulants is that it requires the computation of higher-order integrals of ~the functional under consideration. This is particularly challenging in settings involving long-range spatial correlations. Recently, an important progress was achieved in \cite{LP23,chaos}, where the authors derived a general formula for the cumulants of U-statistics based on Poisson processes, see e.g.~\cite[Theorem~3.2]{chaos} therein. This formula relies on the computation of higher-order integrals, which are particularly accessible for random geometric graphs, where the corresponding integrands are  bounded and of bounded support. This is no longer true in the WRCM where the range of the support depends on the profile function and the weight distribution. 

%
%
While Poisson processes are certainly the simplest model for a random point pattern in Euclidean space, they are inappropriate in settings with substantial spatial correlations. In order to overcome this limitation, one considers determinantal and permanental point processes that are based on a spatial kernel function \cite{hough}. These classes of point processes are particularly well-suited for cumulant-based methods since their cumulant density can be expressed in closed form in terms of the kernel function \cite{brillinger}. One of the key observations in \cite{klein,jolivet} is that the cumulants of U-statistics can again be bounded in terms of higher-order integrals. These integrals are reminiscent to those appearing in \cite[Theorem 3.2]{chaos} but are more complex due to the presence of the kernel function. We show that the results in the Poisson case can be extended to this more complex setting at the cost of weaker bounds on the convergence rate. In this context, we provide a lower bound on the variance of the functional under consideration, which may be of independent interest.

The main contributions of this paper are:
\been
\im[i.] Show asymptotic normality for power-weighted edge lengths and subgraph counts in the WRCM. One of the key achievements of our result is that our bounds are completely transparent in the moment assumptions on the weight distribution and the profile function.
\im[ii.] While \cite{chaos} focusses on the case of Poisson processes, we extend the results to $\alpha$-determinantal processes. This relies on expansions derived in \cite{brillinger,klein,jolivet}. 
\im[iii.] Explain how variance lower bounds can be established for the considered functionals.
\enen

The rest of the manuscript is organized as follows. In Section \ref{sec:mod} we introduce the model and state our main results. Then we discuss certain preliminary results concerning general cumulant computations in Section~\ref{sec:prel}. Sections~\ref{sec:simplex}, \ref{sec:edge} and~\ref{sec:det} contain the proofs of the main results. Finally, Section~\ref{sec:var} contains some considerations regarding the lower variance bounds. 

%
%
\section{Model and main results}
\label{sec:mod}
The main goal of this work is to develop the method of cumulants in the setting of WRCM. One of the core strengths of the cumulant method is that it gives rise to a variety of different limit results as elaborated in \cite{doring2021}. To make our presentation self contained, we reproduce here these consequences.

Let $(X_n)_{n \ge 1}$ be a sequence of square-integrable random variables, let $\g \ge0$ be a constant and let $(\De_n)_{n \ge 1}$ be a positive real-valued sequence. We say that $(X_n)_{n \ge 1}$ satisfies
\begin{itemize}
	\item \textbf{MDP}$(\g, (\De_n)_{n \ge 1})$ if for any positive real-valued sequence $(a_n)_{n \ge 1}$ with
	$$
		\lim_{n \to \ff} a_n=\ff \quad \text{and} \quad \lim_{n \to \ff} \f{a_n}{\De_n^{1/(1+2\g)}}=0
	$$
	the re-scaled random variables $(\widehat X_n)_{n \ge 1} := (a_n^{-1}(X_n-\E[X_n])/\sqrt{\Var(X_n)})_{n \ge 1}$ satisfy a \emph{moderate deviation principle (MDP)} with speed $a_n^2$ and good rate function $\mathcal I(z)=z^2/2$. That is, for all Borel sets $A$,
	$$
		- \inf_{z \in \text{int}(A)} \mathcal I(z) \le \liminf_{n \to \ff} a_n^{-2} \log \P(\widehat X_n \in A) \le \limsup_{n \to \ff} a_n^{-2} \log \P(\widehat X_n \in A) \le -\inf_{z \in \text{cl}(A)} \mathcal I(z),
	$$
	where $\text{int}(A)$ and $\text{cl}(A)$ stand for the interior and the closure of $A$, respectively \cite{ld}. Sometimes this is also referred to as large deviation principle. However, since in our setting the speed $a_n$ is of order smaller than $\sqrt{\ms{Var}(X_n)}$, the term MDP is more appropriate.
	\item \textbf{CI}$(\g, (\De_n)_{n \ge 1})$ if the Bernstein-type concentration inequality
	$$
\P\big(\, |X_n-\E[X_n]|\ge z \sqrt{\Var(X_n)} \,\big) \le 2 \exp \Big( -\f 14 \min\Big\{ \f{z^2}{2^{1+\g}},(z\De_n)^{1/(1+\g)} \Big\} \Big)
	$$
	holds for all $n \ge 1$ and $z \ge 0$,
	\item \textbf{NACC}$(\g, (\De_n)_{n \ge 1})$ if a normal approximation bound with Cr\'amer correction holds, that is, if there exist constants $c_0,c_1,c_2>0$ only depending on $\g$ such that for all $n \ge 1$ and $z \in [0,c_0\De_n^{1/(1+2\g)}]$,
	$$
	 \P\big(\, |X_n-\E[X_n]|\ge z \sqrt{\Var(X_n)} \,\big) = e^{L_{n,z}^+}(1-\Phi(z))\Big(1+c_1\theta_{n,z}^+ \f{1+z}{\De_n^{1/(1+2\g)}}\Big) 
	$$	
	and 
	$$
 \P\big(\, |X_n-\E[X_n]|\le -z \sqrt{\Var(X_n)} \,\big)  = e^{L_{n,z}^-}(1-\Phi(z))\Big(1+c_1\theta_{n,z}^- \f{1+z}{\De_n^{1/(1+2\g)}}\Big) 
	$$
with $\theta_{n,z}^+, \theta_{n,z}^-\in [-1,1]$ and $L_{n,z}^+,L_{n,z}^-\in (-c_2z^3/\De_n^{1/(1+2\g)},c_2z^3/\De_n^{1/(1+2\g)})$, where $\Phi$ is the distribution function of a standard normal random variable.
\end{itemize}  

We now provide a formal definition of the WRCM considered in this paper, confer \cite{glm2,komjathy2,komjathy}.  Let $\PP_n := \{P_i\}_{i \ge 1} := \{(X_i, U_i)\}_{i \ge 1}$ be a Poisson point process on the space $W \ti [0,\ff)$ with intensity measure $t_n \d x \otimes \P_U(\d u)$ for some compact set $W \su \R^d$ with interior points and $t_n >0$ and $\P_U$ is a probability distribution. Then, in the WRCM, between any pair of points $P_i, P_j$ we place an edge independently with probability 
\begin{align}
	\label{eq:prob}
	\vp_n(X_i - X_j, U_i, U_j):=	\vp\bigg(\frac{|B_{\|X_i-X_j\|}|}{\nu_n \k(U_i,U_j)}\bigg), 
\end{align}
where we now explain in detail the individual components. 
\been
\im[i.] $\vp\co [0,\ff) \to [0, 1]$ is the \emph{profile function}, which is non-increasing and for which we assume the normalization $\int \vp(t) \d t = 1$. 
\im[ii.] $|B_{\|X_i-X_j\|}|$ is the volume of the Euclidean ball of radius $\|X_i-X_j\|$ in $\R^d$. Sometimes in literature, the scaling $\|X_i - X_j\|^d$ is used. Since the difference is only in a constant multiple, the two parameterizations are equivalent. For our analysis, the parameterization with the ball-based scaling will simplify substantially the statements and proofs of our main results.
\im[iii.] $(\nu_n)_{n \ge 1}$ is a positive real-valued sequence that is bounded from above. It can be thought of as the volume of a ball around a typical node within which most neighbors of this node are located.
\im[iv.] Finally, $\k\co [0,\ff) \ti [0,\ff) \to [0,\ff)$ is the \emph{interaction kernel}, which is non-decreasing in both arguments. We assume that $\k$ is such that $1\le\k(u,v) \le uv$ for all $u,v \ge 0$. 
\enen

%
%
In a broad sense, our cumulant bounds are similar to those considered earlier in the setting of standard  random geometric graphs \cite{chaos}. However, in the setting of WRCM the convergence will depend sensitively on moment estimates for both the profile function and the weight distribution. To make this precise, we now introduce specific conditions on these quantities. 

We assume that $ U \ge 1 $ a.s.\ and that there are constants $\cuo>0$ and  $\cut\le 1$ such that for all $x\ge1$,
\begin{align}\label{cond:MU}
	M_U(x) :=\mathbb E[U^{x}] \le \cuo^{x} \Ga(1+x)^{\cut}, \tag{\bf MU}
\end{align}
where $\Ga(y):=\int_0^\infty t^{y-1} e^{-t} \d t$ is the standard Gamma function. For instance, if $X$ is exponentially distributed with some parameter $\la>0$, then $U:=X+1$ satisfies \eqref{cond:MU} with $\cut=1$. Moreover, if $X$ is a standard normal random variable and $U \sim X \mid \{X>1\}$, then \eqref{cond:MU} holds with $\cut=1/2$. As we will see in Theorem \ref{thm:simplex} below, the value of the constant $\cut$ is much more important than that of $\cuo$.

%
%
The first functional of interest is the subgraph count, which was also considered in \cite{chaos} and we recall the definition here.  Let $G$ be a fixed connected graph on the vertices $\{1, \dots, q\}$. Then, we let $S_n(G)$ denote $1/q!$ times the number of $q$-tuples $(P_1, \dots, P_q) \in \PP_n^q$ such that for every edge $\{i, j\}$ in $G$, the WRCM has an edge between $P_i$ and $P_j$.

%
%
\bet[Cumulant bounds for subgraph count]
\label{thm:simplex}
Let $ d\ge 1$  and $G$ be a fixed connected graph on $q$ vertices.
Assume that $U \ge 1$ a.s.\ and that condition \eqref{cond:MU} holds and set $A_{\ms{SG}} := (q-1)(1+\cut)$. Moreover, assume that there exists $v > 0$ such that for all $n \ge1$, 
\begin{align}
	\label{eq:pvbg}
	\Var(S_n(G))\ge v|W|t_n^q  \nu_n^{q-1} (1 \vee t_n \nu_n)^{q-1}. \tag{\bf VSG}
\end{align}
Then, $S_n(G)$ satisfies $\emph{\textbf{MDP}}(A_{\ms{SG}},(\b_n)_n)$, $\emph{\textbf{CI}}(A_{\ms{SG}},(\b_n)_n)$ and $\emph{\textbf{NACC}}(A_{\ms{SG}},(\b_n)_n)$ with 
$$\b_n:=\f{\sqrt{v|W|t_n (1 \wedge t_n \nu_n)^{q-1}}}{q^{3q} (b_{\ms{SG}} \vee b_{\ms{SG}}^3/v)}
, \quad\text{where} \quad b_{\ms{SG}}:=((1\vee \cuo^2)(q-1)^{\cut})^{q-1},\quad n \ge 1,$$
\ent

Second, for fixed $\a \ge 0$, we consider the $\a$-power-weighted edge length given by
$$S_n^{(\a)} := \f 12\sum_{P_i\sim P_j \in \PP_n}\| X_i-X_j\|^\a,$$
where $P_i \sim P_j$ indicates the existence of an edge in the WRCM. Here, we additionally need to impose similar assumptions concerning $\vp$ as we did for the weights in \eqref{cond:MU}. More precisely, we assume that there are constants $\cpo>0$ and  $\cpt\le 1$ such that for all $x\ge0$,
\begin{align}
\label{cond:MPHI}
	M_{\vp}(x):= \int_0^\ff \vp(t) t^{x } \, \d t \le  \cpo^{x} \Ga(1+ x)^{\cpt} \quad\text{and}\quad M_{\vp}'(x):=\sup_{t>0} \vp(t) t^{x}\le \cpo^x \Ga(1+ x)^{\cpt}.	 \tag{\bf MPHI}
\end{align}
For instance, \eqref{cond:MPHI} holds if $\vp$ has compact support. Moreover, similarly as for the case of weights, one can also consider the cases of exponential tails where $\cpt = 1$ and Gaussian tails where $\cpt = 1/2$.

%
%
\bet[Cumulant bounds for power-weighted edge lengths]
\label{thm:edge}
Let $ d\ge 1$  and $\a \ge 0$.
Assume that $U \ge 1$ a.s.\ and that conditions \eqref{cond:MU} and \eqref{cond:MPHI} hold and set $A_{\ms{EL}} := 1 + 2\cut(\a_d + 1) + \cpt\a_d$, where $\a_d := \a/d$. Moreover, assume that there exists $v > 0$ such that for all $n \ge1$,  
\begin{align}
	\label{eq:pvbg2}
        \Var(S_n^{(\a)})\ge v|W|t_n^2  \nu_n^{2\a_d + 1} (1 \vee t_n \nu_n). 
\tag{\bf VEL}
\end{align}
Then, the power-weighted edge length $S_n^{(\a)}$ satisfies $\emph{\textbf{MDP}}(A_{\ms{EL}},(\b_n)_n)$, $\emph{\textbf{CI}}(A_{\ms{EL}},(\b_n)_n)$ and $\emph{\textbf{NACC}}(A_{\ms{EL}},(\b_n)_n)$ with 
$$
\b_n:=\f{\sqrt{v|W|t_n (1 \wedge t_n \nu_n)}}{64 (b_{\ms{EL}} \vee b_{\ms{EL}}^3/v)}
,\quad\text{where} \quad b_{\ms{EL}}:=(\cpo (\a_d\vee 1)^{\cpt})^{\a_d}(\cuo (\a_d + 1)^{\cut})^{2(\a_d + 1)},\quad n \ge 1.
$$
\ent

Note that for $U=1$ a.s.~and $\vp(\cdot)=\one_{[0,1]}$, the weight-dependent random connection model reduces to the classical random connection model. In this situation, \cite[Section 4.2]{chaos} establishes the analogue of Theorem \ref{thm:simplex}. Since we can choose $a=\cuo=1$, $\cut=0$ in this case, Theorem \ref{thm:simplex} generalizes \cite[Corollary 4.2]{chaos} to the weight-dependent situation.
	
We also note that for many real-world datasets one can encounter the setting where the weights only have a finite number of moments. In such situations our results do not apply since the condition \eqref{cond:MU} is violated.

Finally in Theorem \ref{thm:det} and \ref{thm:det2} below, we replace the Poisson point process by an $\alpha$-determinantal point process ($\alpha$-DPP) as input. We note that while a homogeneous Poisson point process is a special case of an $\alpha$-DPP, the convergence rates derived for $\alpha$-DPPs are less sharp than the one for Poisson input. Hence, the results in Theorem \ref{thm:det} and \ref{thm:det2} do not supersede the ones in Theorem \ref{thm:simplex} and \ref{thm:edge}.

Recall that for $\alpha\in \R$, a simple point process $\PP$ is called {\em $\alpha$-determinantal point process ($\alpha$-DPP)} with {\em covariance kernel $K:\R^d \times \R^d \to \mathbb C$}, if for all $k \ge 1$, its $k$th product density $\rho_k$ is given by 
\begin{align}
	\r_k(x_1,\dots,x_k)=\sum_{\pi \in \text{Per}(k)} \a^{k-n(\pi)} \prod_{i=1}^k K(x_i,x_{\pi(i)}) \label{eq:rho},
\end{align}
for all $x_1,\dots,x_k \in  \R^d$, where $\text{Per}(k)$ is the set of all permutations of $\{1,\dots,k\}$ and $n(\pi)$ is the number of cycles in the permutation $\pi$.  Important special cases are {\em determinantal processes} ($\a=-1$) and {\em permanental processes} ($\a=1$).
Following \cite[Section 3]{brillinger}, we assume throughout that the kernel is stationary, i.e. $K(x,y)=K(o,y-x)=K_0(y-x)$ for all $x,y \in \R^d$, where $K_0(x)$ is a continuous complex-valued non-negative-definite function, i.e.\ $\sum_{i,j=1}^n z_i K_0(x_i-x_j) \overline{z_j}\ge 0$ for all $x_1,\dots,x_n \in \R^d,\,z_1,\dots,z_n \in \mathbb C$. This implies that $K_0(o)\ge 0$, $K_0(-x)=\overline{K_0(x)}$, $|K_0(x)|\le K_0(o)$ for all $x \in \R^d$ and that the $\alpha$-DPP with kernel $K$ is stationary. Moreover, we assume that
\begin{align}
	\label{eq:DPP}
	K_0(o)\in [0, 1] \quad\text{and}\quad \|K_0\|_1:=\int_{\R^d} |K_0(x)| \,\d x<\infty.
	\tag{\bf DPP}
\end{align}

Under these assumptions, it is shown in \cite[Theorem 2]{brillinger} that the $\alpha$-DPP with kernel $K$ is Brillinger-mixing, which says that for all $k \ge 1$, its $k$th cumulant density is absolutely integrable. 

Since there is no natural analogue to an increasing intensity parameter for $\alpha$-DPPs, we will instead consider the scenario of a stationary point process $\PP$ restricted to an increasing sampling window $W_n:=n^{1/d}W$. This is as in \cite{Yogesh15}. Note that due to the scaling invariance property of a Poisson process, Theorems \ref{thm:simplex} and \ref{thm:edge} could also be formulated in this setup.

%
%
\bet[Cumulant bounds for subgraph count for $\alpha$-DPPs]
\label{thm:det}
Let $ d\ge 1$ and let $\PP$ be an $\alpha$-DPP for some $\alpha \in [-1,1]$. Assume that the  stationary covariance kernel $K$  satisfies the conditions \eqref{eq:DPP}.  Let $G$ be a fixed connected graph with $q$ vertices and consider the subgraph count $S_n(G)$.
 Assume that $U \ge 1$ a.s.\ and that condition \eqref{cond:MU} holds and set $A_{\ms{SG}'} := 2q+(q-1)\cut-1$. Moreover, assume that there exists $v > 0$ such that for all $n \ge1$, 
\begin{align}
	\label{eq:pvbgdet}
	\Var(S_n(G))\ge v |W_n|  \nu_n^{q-1} (1 \vee \nu_n^{q-1}). \tag{\bf VSG'}
\end{align}
Then, the subgraph count $S_n(G)$ satisfies $\emph{\textbf{MDP}}(A_{\ms{SG}'},(\b_n)_n)$, $\emph{\textbf{CI}}(A_{\ms{SG}'},(\b_n)_n)$ and $\emph{\textbf{NACC}}(A_{\ms{SG}'},(\b_n)_n)$ with
$$\b_n:=\f{\sqrt{v|W_n| (1 \wedge  \nu_n)^{q-1}}}{(2q^2)^{3q} (b_{\ms{SG}'} \vee b_{\ms{SG}'}^3/v)}
,\quad\text{where} \quad b_{\ms{SG}'}:=(\|K_0\|_1\vee 1)^{q-1}b_{\ms{SG}},\quad n \ge 1,$$
with $b_{\ms{SG}}$ from Theorem \ref{thm:simplex}.
\ent

The next result is the analogue of Theorem \ref{thm:edge} for $\alpha$-determinantal input.

\bet[Cumulant bounds for power-weighted edge lengths for $\alpha$-DPPs]
\label{thm:det2}
Let $ d\ge 1$ and let $\PP$ be an $\alpha$-DPP for some $\alpha \in [-1,1]$. Assume that the  stationary covariance kernel $K$  satisfies the conditions \eqref{eq:DPP}.  Consider the $\a$-power-weighted edge length $S_n^{(\a)}$.  Assume that $U \ge 1$ a.s.\ and that conditions \eqref{cond:MU} and \eqref{cond:MPHI} hold and set $A_{\ms{EL}'} := 3+\cpt \a_d+2\cut (\a_d+1)$, where $\a_d:=\a/d$. Moreover, assume that there exists $v > 0$ such that for all $n \ge1$, 
\begin{align}
	\label{eq:pvbgdet2}
	\Var(S_n^{(\a)})\ge v|W_n| \nu_n^{2\a_d + 1} (1\vee  \nu_n). 
\tag{\bf VEL'}
\end{align}
Then, the power-weighted edge lengths $S_n^{(\a)}$ satisfies $\emph{\textbf{MDP}}(A_{\ms{EL}'},(\b_n)_n)$, $\emph{\textbf{CI}}(A_{\ms{EL}'},(\b_n)_n)$ and $\emph{\textbf{NACC}}(A_{\ms{EL}'},(\b_n)_n)$ with
$$\b_n:=\f{\sqrt{v|W_n| (1 \wedge  \nu_n)}}{64 (b_{\ms{EL}'} \vee b_{\ms{EL}'}^3/v)}
,\quad\text{where} \quad b_{\ms{EL}'}:=( \|K_0\|_1\vee 1)b_{\ms{EL}},\quad n \ge 1,$$
with $b_{\ms{EL}}$ from Theorem \ref{thm:edge}.
\ent

\begin{remark}
	We stress that the existence of an $\alpha$-DPP with a given covariance kernel $K$ is not guaranteed. Hence, we say that some $\alpha\in \R$ is \emph{admissible}, if there exists a point process whose $k$ product intensity is given by $\eqref{eq:rho}$. By \cite[Theorem 1.2]{shirai}, it is known that all $\alpha\in \{\pm 1/m:\,m \in \N\}$ are admissible.
\end{remark}

\begin{remark}
	 In Propositions \ref{pr:edge} and \ref{pr:simplex} below we show that the variance lower bounds assumed in Theorems \ref{thm:simplex} and \ref{thm:edge} hold for an appropriate sequence $(\nu_n)_{n \ge 1}$  and give a possible choice for the constant $v>0$. Note that variance lower bounds are also required in Theorems \ref{thm:det} and \ref{thm:det2} for $\a \ge 0$. Here, due to the different scenario of an increasing sampling window, we rather work with a sequence where  $(\nu_n/n)_{n \ge 1}$ is bounded. In Remark \ref{rem:det} we address the problem of obtaining variance lower bounds for U-statistics of $\alpha$-determinantal processes if $\alpha <0$.
\end{remark}

%
%
\section{Preliminaries}
\label{sec:prel}
In this section, we set the stage for the proof of our main results and recall some material from \cite{chaos}.  
To keep the presentation self-contained, in Section \ref{sec:partitions}, we start with elementary facts on some partitions which play a fundamental role in the method of cumulants. Then, in Section \ref{sec:cum}, we recall the method of cumulants and following this, the main results from \cite{chaos}. 

%
%
\subsection{Partitions} 
\label{sec:partitions}
Partitions play a central role in our application of the method of cumulants, and in the present section, we reproduce some material from \cite[Chapter 12.2]{last_penrose_lectures}, see also \cite[Section 2.4]{chaos}. However, since we only deal with single U-statistics and not linear combinations of U-statistics, we can simplify the presentation.

Let $m \ge 1$ and let $q \ge 1$. We define $N_0:=0$, $N_\ell:=\ell q $, $\ell \in \{1,\dots,m\}$, and $N:=N_m$,  and put $J_\ell:=\{N_{\ell-1}+1,\dots,N_\ell\}$, $\ell \in \{1,\dots,m\}$. A partition $\s$ of $[N]:=\{1,\dots,N\}$ is a collection $\{B_1,\dots,B_k\}$ of $1 \le k \le N$ pairwise disjoint non-empty sets, called blocks, such that $B_1 \cup \dots \cup B_k=[N]$.  By $\Pi^m(q)$, we denote the set of all partitions $\s$ such that 
$$\max_{B\in\s}\max_{\ell \le m}|B \cap J_\ell|\le 1,$$
and by $\Pi^{m}_{\ge2}(q)$ the set of all $\s\in\Pi^m(q)$ with $\min_{B \in \s}|B|\ge2$.

Every partition $\s \in \Pi^m(q)$ induces a partition $\s^*$ of $\{1,\dots,m\}$ in the following way: $k,\ell \in \{1,\dots,m\}$ are in the same block of $\s^*$ whenever there is a block $B \in \s$ such that $|B \cap J_k|=1$ and $|B \cap J_\ell|=1$. Let $\widetilde \Pi^m(q)$ be the set of all partitions $\s \in \Pi^m(q)$ such that $|\s^*|=1$ and define $\widetilde \Pi^{m}_{\ge2}(q)$ in the same way as $\Pi^m_{\ge2}(q)$.

Let $(\mathbb X,\mathcal X)$ be a measurable space. For functions $f^{(\ell)}:\mathbb X^{q_{\ell}}\to \R$, $\ell \in \{1,\dots,m\}$, we define their tensor product $\otimes_{\ell =1}^m f^{(\ell)}:\mathbb X^N\to \R$ by
$$
(\otimes_{\ell=1}^m f^{(\ell)})(x_1,\dots,x_N):=\prod_{\ell=1}^m f^{(\ell)}(x_{N_{\ell-1}+1},\dots,x_{N_\ell}).
$$
For $\s \in \Pi^m(q)$ the function $(\otimes_{\ell=1} ^{m} f^{(\ell)})_{\s}:\mathbb X^{|\s|}\to \R$ is obtained by replacing in $(\otimes_{\ell=1}^{m} f^{(\ell)})$ all variables that belong to the same block of $\s$ by a new common variable. Note that this way $(\otimes_{\ell=1}^{m} f^{(\ell)})_{\s}$ is only defined up to permutations of its arguments. Since in what follows we always integrate with respect to all arguments of $(\otimes_{\ell=1}^{m} f^{(\ell)})_{\s}$, this does not cause problems.

%
%
\subsection{Method of cumulants} 
\label{sec:cum}
In this section, we introduce the key tool for the proof of the main results, namely the method of cumulants, which originates from \cite{saulis1991}, but also see the survey \cite{doring2021}. Henceforth, $\cum(X_1,\,\dots,\,X_d)$ denotes the joint cumulant of the random variables $(X_1,\,\dots,\,X_d)$. We recall 
that with $i$ the imaginary unit, this is given as 
 $$
 \cum(X_1,\,\dots,\,X_d) := (-i)^m \f{\pa^m}{\pa t_1\cdots \pa t_m}\log \phi_{X_1, \dots, X_m}(t_1, \dots, t_m)\mid_{t_1 = \ldots= t_m = 0},
 $$
 where 
 $$
 \phi_{X_1, \dots, X_m}(t_1, \dots, t_m):=  \E\Big[ \exp\Big( i\sum_{i\le m} t_iX_i \Big) \Big] 
 $$ 
 is the characteristic function of the random vector $(X_1, \dots, X_m)$. We henceforth always assume that all considered random variables have finite moments of all orders.

The following proposition can be regarded as the application of the method of cumulants to Poisson-based U-statistics and is a variant of \cite[Theorem 3.2]{chaos}. In contrast to the latter, we allow for the additional factor $(m!)^p$ in assumption \eqref{cond:ass1} below. Although we do not need it for our purposes, we note that the statement could be generalized to linear combinations of U-statistics.

\bepr[Variant of Theorem 3.2 in \cite{chaos}]
\label{prop:Ustat}
Let $(\PP_n)_n$ be a family of Poisson processes over $\sigma$-finite measure spaces $((\mathbb X_n,\mathcal X_n,\mu_n))_n$ and let $f_n:\mathbb X_n^{q}\to \mathbb R$, $n \in \mathbb N$, be measurable and satisfy the integrability conditions
$$
\max_{i \in \{0,\dots,q\}}\int_{\mathbb X^{i}}\Big(\int_{\mathbb X^{q-i}}f(y_1,\dots,y_i,x_1,\dots,x_{q-i})\, \mu_n^{q-i}(\d(x_1,\dots,x_{q-i}))\Big)^2 \, \mu_n^{i}(\d(y_1,\dots,y_i)) <\ff,
$$
and $\sup_{n \ge 1} \|f_n\|_{L^2(\mu_n^{q})}^2>0$. Define for $n \ge 1$ the random variables
$$
Z_n:= \sum_{(x_1,\dots,x_{q}) \in \PP_{n,\ne}^{q}} f_n(x_1,\dots,x_{q}).
$$
Assume that there is a positive real-valued sequence $(\beta_n)_{n\ge1}$ such that for any $n \in \mathbb N$, for all $\sigma \in \widetilde \Pi_m(q)$ and $m \ge 3$, 
\begin{align}
	\label{cond:ass1}
	{\Var(Z_n)^{-m/2}} \Big|\int_{\mathbb X^{|\s|}} (\otimes_{\ell=1}^m f_n^{(\ell)})_{\sigma}  \, \d\mu^{|\sigma|} \Big|\le (m!)^p \beta_n^{-(m-2)}.
\end{align}
Then $(Z_n)_n$ satisfies \emph{\textbf{MDP}}$(p+q-1,(c_{q}\beta_n^{})_n)$, \emph{\textbf{CI}}$(p+q-1,(c_{q}\beta_n^{})_n)$ and \emph{\textbf{NACC}}$(p+q-1,(c_{q}\beta_n^{})_n)$ with $c_{q}:=1/q^{3q}$.
\enpr

\begin{proof}
	According to \cite[Section~2]{doring2021}, we have to find parameters $\gamma\ge0$ and $\Delta_n>0$, such that the condition
	$$
	\Big|\cum_m\Big(\frac{Z_n}{\sqrt{\Var(Z_n)}}\Big)\Big| \le \frac{(m!)^{1+\gamma}}{\Delta_n^{m-2}}
	$$
	holds, which then implies satisfaction of \textbf{MDP}$(\g, (\De_n)_{n \ge 1})$, \textbf{CI}$(\g, (\De_n)_{n \ge 1})$ and \textbf{NACC}$(\g, (\De_n)_{n \ge 1})$. This condition is widely known as the Statulevi\v{c}ius condition. To do so, we follow the arguments in the proof of \cite[Theorem 3.1]{chaos}. In particular, we rely on the Wiener It\^o integrals $I_q(f)$ as defined in Section~2.1 ibidem.
	By definition of the cumulant, its multilinearity and the triangle inequality we have that
	\begin{align*}
		\Big|\cum_m\Big(\frac{Z_n}{\sqrt{\Var (Z_n)}}\Big)\Big|&=\frac{1}{\Var (Z_n)^{m/2}} |\cum(Z_n,\dots,Z_n)|
		\le \frac{1}{\Var (Z_n)^{m/2}}  |\cum(I_{q}(f_n),\dots,I_{q}(f_n))|.
	\end{align*}
	Now it follows from \cite[Theorem 3.6]{chaos} that
	\begin{align*}
		|\cum(I_{q}(f_n),\dots,I_{q}(f_n))|&\le |\widetilde \Pi^m_{\ge 2}(q)|  \sup_{\sigma \in \tilde \Pi^m_{\ge 2}(q)} \Big|\int_{\mathbb X^{|\s|}} (\otimes_{\ell=1}^m f_n^{(\ell)})_{\sigma}  \, \d\mu^{|\sigma|} \Big|.
	\end{align*}
Thus, we obtain from the first part of \cite[Proposition 6.1]{chaos} that
	$$
	|\widetilde \Pi^m_{\ge 2}(q)| \le |\Pi^m(q)| \le q^{qm} (m!)^q.
	$$
	Putting this together with Assumption~\eqref{cond:ass1} gives that
	$$
	\Big|\cum_m\Big(\frac{Z_n}{\sqrt{\Var (Z_n)}}\Big)\Big|\le q^{qm} (m!)^{p+q} \beta_n^{m-2} \le (m!)^{p+q} (q^{3q} \beta_n)^{m-2}.
	$$
\end{proof}

%
%
\section{Proof of Theorem \ref{thm:simplex}}
\label{sec:simplex}
In this section, we prove Theorem \ref{thm:simplex}, i.e., we derive the asserted bounds for integrals over functions of the form $(f_n^{\otimes m})_\s$ with respect to $\mu_n^{|\s|}$, where for $x_1,\dots,x_q \in (\R^d)^q$, we let $f_n(x_1,\dots,x_q)$ be $1/q!$ times the expected number of subgraphs in the WRCM built on $x_1,\dots,x_q$ that are isomorphic to $G$. While the basic idea is similar to \cite[(7.1)]{chaos}, the presence of long-range correlations through the profile function $\vp$ and the weights requires additional work.

Before bounding the integrals, we need to address a technical point in order to express the considered functional as a U-statistic. This is because in the current framework even when all Poisson points $\{(X_i, U_i)\}_i$ are assumed as given, the decision whether to put an edge between $(X_i, U_i)$ and $(X_j, U_j)$ still involves additional randomness. We resolve this technicality in the same way as in \cite{stein}. More precisely, we formally enlarge the state space of $\PP_n$ from $W \ti [1, \ff)$ to $\mathbb X:=W \ti [1, \ff) \ti [0, 1]^{\mathbb N}$ and augment each Poisson point $(X_i, U_i)$ with additional iid uniform marks $(T_i^{(j)})_{j \ge 1}$. Then, index the points of $\PP_n$ in the lexicographic order and put an edge  between $P_i$ and $P_j$ with $i < j$ if and only if 
$$T_i^{(j)} \le \vp_n(X_i - X_j, U_i, U_j).$$

In the proof, it is useful to note that the partition $\s\in\widetilde{\Pi}^m$ induces a connected graph $(V, E)$ on the vertex set $V=\{1,\dots,|\s|\}$. More precisely, because of the identifications induced by the partition $\s$, the function $(f_n^{\otimes m})_\s$ is a product of certain edge indicators on a collection of $|\s|$ variables. Now, the edge set $E$ consists of those index pairs $i, j \in \{1,\dots,|\s|\}$ appearing in at least one of these edge indicators. Hence, we can write
\begin{align}
\label{eq:int_graphs}
	\int_{\mathbb X^{|\s|}}  (f_n^{ \otimes m})_\s \d\mu_n^{|\s|} = t_n^{|\s|} \int_{W^{|\s|}} \E\Big[\prod_{(i,j)\in E}  \vp_n(x_i - x_j, U_i, U_j)\Big] \d(x_1,\dots,x_{|\s|}) .
\end{align}
Next, we fix a spanning tree of $(V,E)$ with induced edge set $E_s\su E$. Then, we successively integrate the tree $E_s$, where we use spherical coordinates in each integration step. Now, since $\int_0^\ff \vp(t) \d t = 1$, the right-hand side in \eqref{eq:int_graphs} is bounded by
\begin{align}
\label{eq:proof_simplex_integral_step_1}
t_n^{|\s|} v_n^{|\s|-1}|W| \, \E\Big[ \prod_{(i,j)\in E_s} \kappa(U_i,U_j)  \Big].
\end{align}

Let $\operatorname{deg}_{\operatorname{s}}(i)$ denote the degree of $i \in \{1,\dots,|\s|\}$ in $(V,E_s)$. Then, $\k(u,v)\le u  v$ and \eqref{cond:MU} give that 
\begin{align}
	\label{eq:proof_simplex_integral_step_2}
	\E\Big[ \prod_{(i,j)\in E_s} \kappa(U_i,U_j)  \Big] \le   \prod_{i\le|\s|} M_U(\operatorname{deg}_{\operatorname{s}}(i))\le \prod_{i \le |\s|} \cuo^{\operatorname{deg}_{\operatorname{s}}(i)} ( \operatorname{deg}_{\operatorname{s}}(i)!)^{\cut} \;=\; \cuo^{2(|\s|-1)} \prod_{i\le |\s|} \big(\operatorname{deg}_{\operatorname{s}}(i) ! \big)^{\cut}.
\end{align}
Here, the product over $i=1,\dots,|\s|$ can be expressed as
\begin{align*}
	\prod_{i\le|\s|} \big(\operatorname{deg}_{\operatorname{s}}(i) ! \big)^{\cut} = \exp\Big(\cut \sum_{i\le |\s|}\log\Gamma(\operatorname{deg}_{\operatorname{s}}(i)+1)\Big) . 
\end{align*}
First, we show that the function $f(s_1,\dots,s_{|\s|}):=\sum_{i\le |\s|}\log\Gamma(s_i+1)$ is convex. Indeed, for $n\ge1$, the $(n+1)$th  order derivative of $\log\Gamma(s)$ is the polygamma function of order $n$ given by
\begin{align*}
	\Psi_n(s) := (-1)^{n+1} n! \sum_{k\ge 0} \frac{1}{(s+k)^{n+1}} \,.
\end{align*}
This implies that
\begin{align*}
	\frac{\d^2}{\d s^2}\log\Gamma(s)=\Psi_1(s)=\sum_{k\ge 0}\frac{1}{(s+k)^2} >0, 
\end{align*}
which yields convexity of the function $\log\Gamma(s)$. Hence, on the set $\{(s_1,\dots,s_{|\s|}) \in \N^{|\s|}:\,s_1+\cdots+s_{|\s|}=2(|\s|-1)\}$, the function $f$ takes its maximum at an extremal point. In other words, there is some $i_0 \in \{1,\dots,|\s|\}$ such that $s_{i_0} = |\s| - 1$, and $s_i = 1$ for $i \ne i_0$. 
In particular,
$$
	\prod_{i\le|\s|} \big(\operatorname{deg}_{\operatorname{s}}(i) ! \big)^{\cut}\le ((|\s|-1)!)^{\cut}\le ((m(q-1))!)^{\cut} \le (q-1)^{m(q-1)\cut} (m!)^{(q-1)\cut},
$$
where we have used $|\s| \le m(q-1)+1$ and the elementary inequality 
\begin{align}
	(ck)! \le c^{ck} (k!)^c,	\label{eq:facest}
\end{align}
which holds for all $c\ge 1$ and $k \in \N$. Thus,
\begin{align}
	\label{eq:kap}
	\E\Big[ \prod_{(i,j)\in E_s} \kappa(U_i,U_j)  \Big] \le \cuo^{2(|\s|-1)}(q-1)^{m(q-1)\cut} (m!)^{(q-1)\cut} \le b_{\ms{SG}}^m (m!)^{(q-1)\cut},
\end{align}
where we recall  that $b_{\ms{SG}}:=((1\vee \cuo^2)(q-1)^{\cut})^{q-1}$ from Theorem \ref{thm:simplex}. Hence, we arrive at
	\begin{align*}
		\int_{\mathbb X^{|\s|}}  (f_n^{ \otimes m})_\s \d\mu_n^{|\s|} \le t_n^{|\s|} v_n^{|\s|-1}|W| \, b_{\ms{SG}}^m(m!)^{(q-1)\cut}.
	\end{align*}

Now, by Condition \eqref{eq:pvbg}, there exists some constant $v>0$ such that
$$
\Var(S_n(G)) \ge v|W|t_n^q \nu_n^{q-1} (1 \vee t_n^{q-1}\nu_n^{q-1}).
$$
Next, we distinguish the cases $ t_n \nu_n \le 1$ and $ t_n \nu_n>1$. 
If $t_n\nu_n \le 1$,  we use that $|\s|\ge q$ and find with $a:=(q-1)\cut$ that
\begin{align*}
	&(m!)^{-a}\big(\Var S_n(G)\big)^{-m/2} 	\int_{\mathbb X^{|\s|}}  (f_n^{\otimes m})_\s \,\d\mu_n^{|\s|} \\
	&\le \frac{|W| t_n^q b_{\ms{SG}}^m \nu_n^{q-1}}{|W|^{m/2}v^{m/2}t_n^{mq/2}\nu_n^{m(q-1)/2}} =\frac{b_{\ms{SG}}^m |W|^{-(m-2)/2}}{v^{m/2} t_n^{(m-2)q/2}  \nu_n^{(m-2)(q-1)/2}}  \le  (1 \vee b_{\ms{SG}}^2/v)^{m-2} \Big(b_{\ms{SG}}/\sqrt{|W| v t_n^q \nu_n^{q-1}}\Big)^{m-2}. 
\end{align*}
If $t_n\nu_n \ge 1$, we use the bound $|\s|\le m(q-1)+1$ and obtain that
\begin{align*}
	&(m!)^{-a}\big(\Var  S_n(G)\big)^{-m/2}	\int_{\mathbb X^{|\s|}}  (f_n^{\otimes m})_\s \,\d\mu_n^{|\s|} \\
		&\le \frac{|W| t_n^{m(q-1)+1} b_{\ms{SG}}^m \nu_n^{m(q-1)}}{|W|^{m/2}v^{m/2}t_n^{mq-m/2}\nu_n^{m(q-1)}}=\frac{|W|^{-(m-2)/2}t_n^{1-m/2}b_{\ms{SG}}^m}{v^{m/2}} \le (1 \vee b_{\ms{SG}}^2/v)^{m-2} (b_{\ms{SG}}/\sqrt{|W| v t_n })^{m-2}.
\end{align*}
Therefore, a cumulant bound as in \eqref{cond:ass1} is satisfied and the result follow from Proposition~\ref{prop:Ustat}. 

%
%
\section{Proof of Theorem \ref{thm:edge}}
\label{sec:edge}

In this section, we prove Theorem \ref{thm:edge}. We apply Proposition \ref{prop:Ustat} with state space $\mathbb X:=W \ti [1, \ff) \ti [0, 1]^{\mathbb N}$, $k=1$, $q=2$ and 
\begin{equation*}
	f_n(P_i, P_j):= |B_{\|X_i - X_j\|}|^\a \one\Big\{T_i^{(j)}\le \vp_n(X_i - X_j, U_i, U_j)\Big\}.
\end{equation*}
   In the proof, we will replace the $\a$-power-weighted edge length $\|X_i-X_j\|^{\a}$ by $|B_{\|X_i-X_j\|}|^{\a}=|B_1|^{\a} \|X_i-X_j\|^{d \a}$.  Note that the linear scaling of $S_n^{(\a)}$ by $|B_1|^{\a}$ has no influence on the fraction $S_n^{(\a)}/\sqrt{\Var(S_n^{\tau})}$ and, therefore, does not affect our bounds. As we will see below, this reparametrization will make the integral computations substantially more accessible.

As in Section \ref{sec:simplex}, we use that the partition $\s$ induces a connected graph $(V, E)$ on the vertex set $V=\{1,\dots,|\s|\}$. In fact, in the present setting this graph is far easier to comprehend since $q = 2$. Moreover, for $(i,j)\in E$ we henceforth  let $m_{ij}$ be the multiplicity of $(i,j)$ in $\s$ and note that $\sum_{(i,j)\in E}m_{ij}=m$. Then,
\begin{align}
\label{eq:int_edges}
	\int_{\mathbb X^{|\s|}}  (f_n^{\otimes m})_\s \d\mu_n^{|\s|} = t_n^{|\s|}\E\Big[\int_{W^{|\s|}} \prod_{(i,j)\in E} |B_{\|x_i-x_j\|}|^{\a m_{ij}}  \vp_n(x_i - x_j, U_i, U_j) \Big] \d(x_1,\dots,x_{|\s|}) .
\end{align}
Again, we fix the edge set $E_s \su E$ of a spanning tree and rewrite the inner integral in the above as
\begin{align*}
	&\int_{W^{|\s|}} \hspace{-0.1cm} \bigg[ \prod\limits_{(i,j)\in E_s} \vp_n(x_i - x_j, U_i, U_j)) |B_{\|x_i-x_j\|}|^{\a m_{ij}} \bigg]   \bigg[ \prod\limits_{(i,j)\in E \sm E_s} \vp_n(x_i - x_j, U_i, U_j)  |B_{\|x_i-x_j\|}|^{\a m_{ij}} \bigg] \d(x_1,\dots,x_{|\s|})\\
	&\quad\le \bigg[\prod\limits_{(i,j)\in E \setminus E_s} (\nu_n \k(U_i,U_j))^{{\a m_{ij}}{}} M_{\vp}'(\a m_{ij})\bigg] \int_{W^{|\s|}} \prod\limits_{(i,j)\in E_s} \vp_n(x_i-x_j,U_i,U_j) |B_{\|x_i-x_j\|}|^{\a m_{ij}}\d(x_1,\dots,x_{|\s|})
\end{align*}
where $M_{\vp}'$ is given at \eqref{cond:MPHI}. Next, we bound the integral over $x_1,\dots,x_{|\s|}$. Again, successively integrating out the contributions of the spanning tree leaves gives the bound
$$
 \int_{W^{|\s|}} \prod\limits_{(i,j)\in E_s} \vp_n(x_i-x_j,U_i,U_j) |B_{\|x_i-x_j\|}|^{\a m_{ij}}\d(x_1,\dots,x_{|\s|})\le |W| \prod_{(i,j) \in E_s} \Big[(\nu_n \k(U_i,U_j))^{{\a m_{ij}}{}+1} M_{\vp}(\a m_{ij})\Big].
$$
Now, Assumption \eqref{cond:MPHI} gives that 
$$\prod\limits_{(i,j)\in E \setminus E_s}M_{\vp}'(\a m_{ij}) \prod\limits_{(i,j)\in E_s}  M_{\vp}(\a m_{ij}) \le \cpo^{\a m} \Ga(1 +  \a m)^{\cpt}.$$
Thus, using the H\"older inequality, $|\s|\le m+1$ and \eqref{eq:facest} shows that \eqref{eq:int_edges} is bounded by
\begin{align*}
	&t_n^{|\s|}|W| \cpo^{{\a m}{}} \Ga(1 + \a m)^{\cpt}\E\Big[  \prod_{(i,j) \in E_s} (\nu_n \k(U_i,U_j))^{{\a m_{ij}}{}+1} \prod_{(i,j) \in E \setminus E_s} (\nu_n \k(U_i,U_j))^{{\a m_{ij}}{}} \Big]\\
	&\le t_n^{|\s|}|W| \cpo^{{\a m}{}} \Ga(1 +\a m)^{\cpt}  \nu_n^{{\a m}{}+|\s|-1}\E\big[  \k(U_i,U_j)^{(\a + 1) m} \big]\\
	&\le t_n^{|\s|}|W| (\cpo (\a \vee 1)^{\cpt})^{{\a m}{}} (m!)^{\cpt\a }  \nu_n^{{\a m}{}+|\s|-1}\E\big[  \k(U_i,U_j)^{(\a + 1) m} \big],
\end{align*}
where we have used that $\Ga(1+\a m) \le (\a \vee 1)^{\a m} (m!)^{\a}$ in the last inequality, which follows from convexity for $\a\le 1$ and from \eqref{eq:facest} for $\a >1$. Hence, it remains to bound the expected value on the right-hand side. To do so, since $\k(u,v) \le uv$, using condition \eqref{cond:ass1} we find with the bound \eqref{eq:facest} that 
\begin{align}\label{eq:Kubou}
	\E\big[  \k(U_i,U_j)^{(\a+1) m} \big]^{1/2} \le \cuo^{(\a+1) m}\Ga(1+ (\a+1) m)^{\cut}&\, \le (\cuo (\a+1)^{\cut})^{(\a+1) m}(m!)^{\cut (\a+1)}.
\end{align}
 Moreover, by Condition  \eqref {eq:pvbg2}, there exists some $v>0$ such that
$$
	\Var(S_n^{(\a d)})\ge v|W|t_n^2  \nu_n^{2\a + 1} (1 \vee t_n \nu_n).
$$
Now we combine these bounds and consider the cases $t_n \nu_n\ge 1$ and $t_n \nu_n <1$ separately. Recall that $b_{\ms{EL}}:=(\cpo (\a\vee 1)^{\cpt})^{{\a }{}}(\cuo (\a + 1)^{\cut})^{2(\a + 1)}$ and $a:=2\cut (\a+1) + \cpt \a$. If $t_n \nu_n <1$, we use that $|\s|\ge 2$ and obtain that  

\begin{align*}
	&(m!)^{-a}\Var(S_n^{(\a d)})^{-m/2} \, 	\int_{\mathbb X^{|\s|}} ( f_n^{\otimes m})_\s \d\mu_n^{|\s|}\\
	 &\quad \le \frac{|W|t_n^2 \nu_n^{\a m} b_{\ms{EL}}^m \nu_n } {|W|^{m/2}v^{m/2} \nu_n^{\a m} t_n^{m} \nu_n^{m/2} }
	=\f{b_{\ms{EL}}^m }{|W|^{(m - 2)/2}v^{m/2} t_n^{m-2}\nu_n^{(m-2)/2}} \le (1 \vee b_{\ms{EL}}^2/v)^{m-2} \Big(b_{\ms{EL}}/\sqrt{|W| v t_n^2 \nu_n}\Big)^{m-2}.
\end{align*}
If $t_n \nu_n \ge1$, because of $|\s|\le m+1$, we have that
 \begin{align*}
	&(m!)^{-a} \Var(S_n^{(\a d)})^{-m/2} \, 	\int_{\mathbb X^{|\s|}} ( f_n^{\otimes m})_\s \d\mu_n^{|\s|}\\
	&\quad  \le \frac{|W|t_n^{m+1} \nu_n^{\a m} b_{\ms{EL}}^m   \nu_n^m } {|W|^{m/2}v^{m/2} \nu_n^{\a m}t_n^{2m-m/2} \nu_n^{m} }
	 =\f{b_{\ms{EL}}^m  }{|W|^{(m - 2)/2}t_n^{(m - 2)/2}v^{m/2}}\le (1 \vee b_{\ms{EL}}^2/v)^{m-2} \Big(b_{\ms{EL}}/\sqrt{|W| v t_n}\Big)^{m-2}.
 \end{align*}
\hfill$\square$

%
%
\section{Proof of Theorems \ref{thm:det} and \ref{thm:det2}}
\label{sec:det}
Before giving the proof of Theorem \ref{thm:det}, we provide a formula from the Appendix of \cite{klein} for the $m$th cumulant of a $U$-statistic of an $\alpha$-DPP and explain how this formula can be interpreted graphically. This will allow us to proceed similarly as in the proof of Theorem \ref{thm:simplex}, with the difference that edges in the graph induced by $\s$ can now have two different meanings.

Given a partition $\s \in \Pi^m(q)$ of $\{1,\dots,mq\}$, let $\r$ be a partition of $\{1,\dots,|\s|\}$. Recalling the definition of $J_t$ from Section \ref{sec:partitions}, we call the pair $(\s,\r)$ {\em decomposable}, if there are a decomposition $\{T_1,T_2\}$ of $\{1,\dots,q\}$ and decompositions $\s=\s^{(1)} \cup \s^{(2)}$ and $\r=\r^{(1)} \cup \r^{(2)}$ such that $\s^{(1)}$ is a partition of $\bigcup_{t \in T_1} J_t$, $\s^{(2)}$ is a partition of $\bigcup_{t \in T_2} J_t$, $\r^{(1)}$ is a partition of $T_1$ and $\r^{(2)}$ is a partition of $T_2$, where we . Otherwise, $(\s,\r)$ is {\em indecomposable}. Given $\s \in  \Pi^m(q)$, let $\Pi(q;\s)$ be the set of all partitions $\r$ such that the pair $(\s, \r)$ is indecomposable. This notion of indecomposability is compatible with the notion of indecomposable integrals used in \cite{klein} in the following sense. An integral is indecomposable in the sense of \cite{klein}, if it cannot be defined with respect to a decomposable partition. Hence, the notion of indecomposable partitions if less restrictive than indecomposable integrals. 

%
%
Therefore, with $f_n$ defined as in the proof of Theorem \ref{thm:simplex}, we have
\begin{align}
	\label{eq:dek}
	|\cum_m(S_n(G))|\le \sum_{\s \in \Pi^m(q)} \sum_{\rho \in \Pi(q;\s)}  \int_{\mathbb X^{|\sigma|}}  (f_n^{\otimes m })_{\s}  \, \prod_{Q \in \r} |c^{(|Q|)}| \,\d \mu_n^{|\s|},
\end{align}
where 
\begin{align}
	\label{eq:det}
	c^{(k)}(x_1, \dots, x_k) = \alpha^{k - 1} \sum_{\pi \in \text{Per}(k)} K({x_1, x_{\pi(2)}})K({x_{\pi(2)}, x_{\pi(3)}})\cdots K({x_{\pi(k)}, x_{1}}),
\end{align}
is the $k$th cumulant density of the $\alpha$-DPP and we recall that $\text{Per}(k)$ is the set of all permutations of $1,\dots,k$.

As above, we note that the choice of partitions $\Pi^m(q)$ and $\Pi(q;\s)$ induces a graph on $\{1, \dots, |\s|\}$. However, in contrast to the Poisson setting, an edge in this graph can now have two different meanings, and therefore we refer to them as type I or type II edges.  First, an edge of type I is put between nodes $i$ and $j$ if the product of the functions $f_n$ imposes an edge. This is just as in Theorem \ref{thm:simplex} and \ref{thm:edge}. However, secondly, an edge of type II is now also if $i$ and $j$ are members of the same element in $\rho \in \Pi(q;\s)$. This is new in the determinantal case. Then, the indecomposability assumption means that the so-constructed graph is connected.

\begin{proof}[Proof of Theorem \ref{thm:det}]
First we bound the number of partitions involved in the double sum in \eqref{eq:dek}. From \cite[Proposition 6.1]{chaos} we have that $|\Pi^m(q)|\le q^{qm}(m!)^q$. To bound the sum over the partitions $\Pi(q;\s)$, we observe from \eqref{eq:det} that when expanding the $k$th order cumulant measure, we obtain a sum consisting of $(m - 1)! \le m!$ summands. Hence, given $\s \in \Pi^m(q)$, the number of terms  \eqref{eq:det} is bounded by 
$$\sum_{r=1}^{|\s|}\sum_{\substack{\rho \in \Pi(q;\s)\\ |\rho|=r}} |Q_1|! \cdots |Q_r|! \le \sum_{r=1}^{|\s|} \f{2^{|\s|-1} |\s|!}{r!} \le  2^{|\s|} |\s|!\le (2q)^{qm} (m!)^q,$$
where the upper bound of the inner sum holds by \cite[Lemma 3.5]{raic} and we have used \eqref{eq:facest} in the last step. Since $|\alpha| \le 1$, this gives that $|\cum(S_n(G))|$ is bounded by
\begin{align}
(2q^2)^{qm} (m!)^{2q} \sup_{\substack{\s \in  \Pi^m(q)\\ \rho \in \Pi(q;\s)}}&\int_{\mathbb X^{|\sigma|}}(f_n^{\otimes m })_{\s}   \prod_{Q \in \r}\sup_{\pi \in \text{Per}(|Q|)} |K({x_1, x_{\pi(2)}})K({x_{\pi(2)}, x_{\pi(3)}})\cdots K({x_{\pi(|Q|)}, x_{1}})|\, \mu_n^{|\s|}(\d \xx). \label{eq:cumbou}
\end{align}
Hence, it remains to bound the individual integrands in \eqref{eq:cumbou}, which consist of a product of the functions $f(\cdot)$ and products of the kernels $K(\cdot, \cdot)$. 

	We construct a spanning tree such that for points $x_i$, $x_j$ in the same element of the partition $\rho$, we place a type II edge $(x_i,x_j)$ in the spanning tree if the factor $K(x_i,x_j)$ is present in the integral. We put type I edges in the spanning tree that connect the different $Q \in \r$. To bound the integral  \eqref{eq:cumbou}, we integrate over the edges of the spanning tree.  If we integrate over a type I edge $((x,u), (y,v))$, then we obtain a contribution of $K_0(o) \nu_n \k(u,v)$. This is analogous to the proof of Theorem \ref{thm:simplex}. Second, each time we integrate over a type II edge, we obtain a contribution bounded by $\|K_0\|_1$. This follows from the integrability of the covariance kernel $K$. Let $E_{s,1}$ denote the set of type I in the spanning tree and $E_{s,2}$ denote the set of type II edges in the spanning tree, respectively. Since $\sup_{x, y \in \R^d}|K(x,y)|\le K_0(o)\le 1$, we deduce from \eqref{eq:kap} the bound
\begin{align*}
&\E \Big[ \int_{W_n^{|\s|}} \prod_{(i,j)\in E_{s,1}} \varphi_n(x_i-x_j,U_i,U_j) \prod_{(i,j)\in E_{s,2}} |K(x_i, x_j)|  \, \d (x_1,\dots,x_{|\s|}) \Big] \\
	&\quad\le |W_n| \nu_n^{|E_{s,1}|} \E\Big[ \prod_{(i,j)\in E_{s,1}} \k(U_i,U_j) \Big] \|K_0\|_1^{|E_{s,2}|} \le |W_n|  \nu_n^{|E_{s,1}|} (m!)^{(q-1)\cut} b_{\ms{SG}}^m \|K_0\|_1^{|E_{s,2}|}.
\end{align*}
where $b_{\ms{SG}}$ is given in Theorem \ref{thm:simplex}. Since $q-1 \le |E_{s,1}|\le m(q-1)$ and $|E_{s,2}|\le m-1$ we can in general conclude that
$$
	|\cum_n(S_n(G))|	\le |W_n| (1 \vee \nu_n)^{m(q-1)}  ((\|K_0\|_1\vee 1)(2q^2)^qb_{\ms{SG}})^m  (m!)^{2q+a}.
$$
From here, the assertion follows as in the proof of Theorem \ref{thm:simplex}.
\end{proof}

\begin{proof}[Proof of Theorem \ref{thm:det2}]
	To prove Theorem \ref{thm:det2}, we combine arguments used in the proofs of Theorems \ref{thm:edge} and \ref{thm:det}. Again, we replace the $\a$-power-weighted edge length $\|X_i-X_j\|^\a$ by $|B_{\|X_i-X_j\|}|^\tau$ to ease the presentation. For $m \ge 1$ let $\s \in \Pi^m(2)$ and $\r \in \Pi(2;\s)$. Let $(V,E)$ be the connected graph induced by the pair $(\s,\r)$ on the vertex set $V=\{1,\dots,|\s|\}$. As above, we can partition $E$ into a set $E_1$ of type I edges and a set $E_2$ of type II edges. For $(i,j)\in E_1$ let $m_{ij}$ be the multiplicity of $(i,j)$ in $(V,E)$ and note that $\sum_{(i,j)\in E_1}m_{ij}=m$. Let $E_{s,2}\su E_2$ be such that the induced subgraph of $(V,E)$ with edge set $E_1 \cup E_{s,2}$ is connected. Since $\sup_{x, y \in\R^d}|K(x,y)| \le K_0(o)\le 1$ for all $x,y \in W_n$, we obtain analogously to \eqref{eq:cumbou} that $|\cum_m(S_n^{(\a d)})|$ is bounded by
\begin{align*}
64^m (m!)^4 \sup_{\s \in  \Pi^m(q)} \sup_{\rho \in \Pi(q;\s)}\E \Big[ \int_{W_n^{|\s|}} \prod_{(i,j)\in E_1} \varphi_n(x_i-x_j,U_i,U_j) |B_{\|x_i-x_j\|}|^{\a} \prod_{(i,j)\in E_{s,2}} |K(x_i, x_j)|  \, \d (x_1,\dots,x_{|\s|}) \Big].
\end{align*}
Now we fix some $E_{s,1} \su E_1$ such that $E_{s,1} \cup E_{s,2}$ is the edge set of a spanning tree of $(V,E)$. Then we find as in the proof of Theorem \ref{thm:edge} that $|\cum_m(S_n^{(\a d)})|$ is bounded by
\begin{align}
64^m (m!)^4	& \E\bigg[\prod\limits_{(i,j)\in E_1 \setminus E_{s,1}} (\nu_n \k(U_i,U_j))^{{\a m_{ij}}{}} M_{\vp}'(\a m_{ij}) \nonumber\\
	&\quad \times  \int_{W_n^{|\s|}} \prod_{(i,j)\in E_{s,1}} \varphi_n(x_i-x_j,U_i,U_j) |B_{\|x_i-x_j\|}|^{\a} \prod_{(i,j)\in E_{s,2}} |K(x_i, x_j)|  \, \d (x_1,\dots,x_{|\s|}) \bigg], \label{eq:int_edges:det2}
\end{align}
where $M_{\vp}'$ is given at \eqref{cond:MPHI}. Next, we bound the integral over $x_1,\dots,x_{|\s|}$. Again, successively integrating out the contributions of the leaves gives that the expectation above is bounded by
$$
|W_n| \|K_0\|_1^{|E_{s,2}|}\prod_{(i,j) \in E_{s,1}} \Big[(\nu_n \k(U_i,U_j))^{{\a m_{ij}}{}+1} M_{\vp}(\a m_{ij})\Big].
$$
Now, Assumption \eqref{cond:MPHI} and $\sum_{(i,j) \in E_1} m_{ij}= m$ give that 
$$\prod\limits_{(i,j)\in E_1 \sm E_{s,1}}M_{\vp}'(\a m_{ij}) \prod\limits_{(i,j)\in E_{s,1}}  M_{\vp}(\a m_{ij}) 
\le \cpo^{\a m} \Ga(1 +  \a m)^{\cpt}.$$
Moreover, by the H\"older inequality,
$$
\E\Big[  \prod_{(i,j) \in E_{s,1}} (\nu_n \k(U_i,U_j))^{{\a m_{ij}}{}+1} \prod_{(i,j) \in E_1 \setminus E_{s,1}} (\nu_n \k(U_i,U_j))^{{\a m_{ij}}{}} \Big] \le \nu_n^{\a m+|E_{s,1}|}\E\big[  \k(U_i,U_j)^{(\a + 1) m} \big].
$$
Thus, \eqref{eq:int_edges:det2} is bounded by
$$
64^m (m!)^4 |W_n| \|K_0\|_1^{|E_{s,2}|} \cpo^{{\a m}{}} \Ga(1 + \a m)^{\cpt}\nu_n^{\a m+|E_{s,1}|}\E\big[  \k(U_i,U_j)^{(\a + 1) m} \big].
$$
 Using that $\E\big[  \k(U_i,U_j)^{(\a+1) m} \big]^{1/2}  \le (\cuo (\a+1)^{\cut})^{(\a+1) m}(m!)^{\cut (\a+1)}$ by \eqref{eq:Kubou} and the bounds $q-1 \le |E_{s,1}|\le m(q-1)$ and $|E_{s,2}|\le m-1$, we arrive at the estimate
$$
\big|\cum_m(S_n^{(\a d)})\big| \le  |W_n|   ( (1 \vee\|K_0\|_1) 64b_{\ms{EL}})^m   \nu_n^{\a m} (1\vee \nu_n)^{m(q-1)}(m!)^{4+\cpt \a+2\cut (\a+1)},
$$ 
where $b_{\ms{EL}}:=(\cpo (\a\vee 1)^{\cpt})^{{\a }{}}(\cuo (\a + 1)^{\cut})^{2(\a + 1)}$. Now, the claim follows as in the proofs above.
\end{proof}

%
%
\section{Variance lower bounds}
\label{sec:var} 
All of our main results depend on lower bounds for the variance of $S_n(G)$ and $S_n^{(\a)}$, respectively. In this final section, we elaborate on how these lower bounds can be established. Although the following propositions are stated for the scenario of a fixed point process restricted to a growing sampling window $W_n$, they analogously hold in the case of a Poisson process with a growing intensity $t_n$ on a fixed sampling window $W$.

For some $a>0$ let
$$W_n^{-a}:=\{x \in W_n:\,B_{(a/|B_1|)^{1/d}}(x) \su W_n\}$$
be the set of all $x \in W_n$ such that a ball centered at $x$ with volume $a$ is still contained entirely in $W_n$. 

We start with the the lower bound on $\Var(S_n^{(\a d)})$. 
%
%
\bepr[Variance lower bounds for power-weighted edge count]
\label{pr:edge}
For $\alpha \in [0,1]$ let  $\PP$ be a stationary $\alpha$-DPP. Let $\e > 0$ and $\a \ge 0$. Then for $n \ge 1$,
$$
\Var(S_n^{(\a d)})\ge v|W_n| \nu_n^{2\a + 1} (1\vee  \nu_n),
$$
where 
$$
v:=\inf_{n \ge 1} \f{|W_n^{-\e\nu_n}|}{|W_n|} |B_1|^{-2\a} K_0(o)^2 \int_0^{\e} t^{2\a + 1} \vp(t) \d t \Big( 1\wedge 2 K_0(o) \int_0^{\e} t^{2\a + 1} \vp(t) \d t \, \Big) .
$$
\enpr

Note that if the sequence $(\nu_n/n)_{n \ge 1}$ is bounded, we have that $v>0$. Indeed, our assumptions (that $\int_0^\infty \vp(t) \d t=1$ and that $\vp$ is non-increasing) guarantee that there exists some $\e>0$ such that $\vp(\e)>0$ and $\inf_{n \ge 1}  \f{|W_n^{-\e \nu_n}|}{|W_n|}>0$.

\bep[Proof of Proposition \ref{pr:edge}]
Expanding the second moment of $S_n^{(\a d)}$ into three terms, we find that  the variance of $S_n^{(\a d)}$ is given by
\begin{align*}
	|B_1|^{2\a}\Var(S_n^{(\a d)}) &=|B_1|^{2\a}\big(\E[(S_n^{(\a d)})^2] - \E[S_n^{(\a d)}]^2 \big)\\
	&= \int_{W_n^2} |B_{\|x-y\|}|^{2\a}\E\big[\vp_n(x - y, U, V)\big] \r_2(x,y) \d(x,y)\\
	&\,+2\int_{W_n^3} |B_{\|x-y\|}|^{\a} |B_{\|x-y'\|}|^{\a} \E\big[\vp_n(x - y, U, V)\vp_n(x - y', U, V')\big] \r_3(x,y,y') \d(x,y,y')\\
	&\,+ \int_{W_n^4} |B_{\|x-y\|}|^{\a} |B_{\|x'-y'\|}|^{\a} \E\big[\vp_n(x - y, U, V)\vp_n(x'-y',U',V')]\r_4(x,x',y,y')) \d (x,y,x',y')\\
	&\,- \int_{W_n^4} |B_{\|x-y\|}|^{\a} |B_{\|x'-y'\|}|^{\a} \E\big[\vp_n(x - y, U, V)\vp_n(x' - y', U', V')\big] \r_2(x,y)\r_2(x',y') \d (x,y,x',y'),
\end{align*}
where $U, U',V, V' \sim \P_U$ are independent. Since for an $\alpha$-DPP with $\alpha \ge0$ all terms on the right-hand side in \eqref{eq:rho} are non-negative, it holds that $\rho_{4}(x,x',y,y')\ge \rho_2(x,y) \rho_2(x',y')$ for all $x,x',y,y' \in \R^{d}$. Therefore,
\begin{align*}
	|B_1|^{2\a}	\Var(S_n^{(\a d)}) &\ge \int_{W_n^2} |B_{\|x-y\|}|^{2\a}\E\big[\vp_n(x - y, U, V)\big] \r_2(x,y)\d(x,y)\\
	&\quad+2\int_{W_n^3} |B_{\|x-y\|}|^{\a} |B_{\|x-y'\|}|^{\a} \E\big[\vp_n(x - y, U, V)\vp_n(x - y', U, V')\big] \r_3(x,y,y')\d(x,y,y').
\end{align*}
We now give lower bounds for the two terms separately, starting with the first one. We use spherical coordinates and exploit that $\k(u,v)\ge 1$. Since $\rho_2(x,y)\ge K_0(o)^2$ for all $x,y \in \R^d$, this gives
\begin{align*}
	\int_{W_n^2} |B_{\|x-y\|}|^{2\a}\E\big[\vp_n(x - y, U, V)\big] \r_2(x,y)\d(x,y) &\ge K_0(o)^2 |W_n^{-\e \nu_n}| \nu_n^{2\a + 1} \int_0^{\e} t^{2\a + 1} \vp(t) \d t.
\end{align*}
Similarly, for the second integral we obtain the lower bound
$$K_0(o)^3|W_n^{-\e \nu_n}| \, \nu_n^{2(\a+1)}\Big(\int_0^{\e} t^{2\a + 1} \vp(t) \d t\Big)^2.
$$
Hence, combining the two bounds yields the assertion.
\enp

Next, we consider the variance of the subgraph count $S_n(G)$.

%
%

\bepr[Variance lower bounds for the subgraph count]
\label{pr:simplex}
For $\alpha \in [0,1]$ let  $\PP$ be a stationary $\alpha$-DPP, let $G$ be a fixed connected graph on $q$ vertices. Let $\varepsilon>0$. Then for $n \ge 1$,
\begin{align}
	\Var(S_n(G))\ge v |W_n| \nu_n^{q-1} (1 \vee \nu_n^{q-1}), \label{eq:pvbg1}
\end{align}
where 
$$
v:=\inf_{n \ge 1}  \f{|W_n^{-\e \nu_n/2^d}|}{|W_n|} \vp(\e)^{2 \binom q2} K_0(o) (\e K_0(o)/2^d)^{q-1} (1 \wedge q(\e K_0(o)/2^d)^{q-1}).
$$
\enpr

\bep
Let $f_n$ be defined as in the proof of Theorem \ref{thm:simplex}. Moreover, for $k \in \N$ let $\rho_\ell$ be the $k$th product density of $\PP$. Then, we have
\begin{align}
	\Var(S_n(G)) =& \sum_{k=0}^{q} \binom qk \int_{\mathbb X^q} \int_{\mathbb X^k} f_n(\mathbf s) f_n(\mathbf s_{q-k},\mathbf t)  \rho_{q+k}(\mathbf x,\mathbf y) \mu_n^k(\d \mathbf t) \mu_n^q(\d \mathbf s) \nonumber\\
	&\quad- \int_{\mathbb X^q} \int_{\mathbb X^q} f_n(\mathbf s) f_n(\mathbf t) \rho_q(\xx)\rho_q(\yy) \mu_n^q(\d \mathbf t) \mu_n^q(\d \mathbf s),\label{eq:vareq}
\end{align}
where $\mathbf s:=(s_1,\dots,s_q)$ and $x_i$ is the first coordinate of $s_i$ for $i=1,\dots,q$, and $\mathbf s_{q-k}:=(s_1,\dots,s_{q-k})$ (similarly for $\mathbf t$ and $\yy$). Since for $\alpha>0$ all terms on the right-hand side in \eqref{eq:rho} are non-negative, it holds that $\rho_{2q}(\xx,\yy)\ge \rho_q(\xx) \rho_q(\yy)$ for all $\xx \in \R^{qd}$, $\yy \in \R^{qd}$. This yields
\begin{align}
	\Var(S_n(G)) \ge \sum_{k=0}^{q-1} \ \binom qk \int_{\mathbb X ^q} \int_{\mathbb X ^k} f_n(\mathbf s) f_n(\mathbf s_{q-k},\mathbf t)  \rho_{q+k}(\xx,\yy) \mu_n^k(\d \mathbf t) \mu_n^q(\d \mathbf s).\label{eq:varineq}
\end{align}
To bound the term with $k=0$, note that $f_n(\xx)$ is bounded by the indicator that $x_2,\dots,x_q$ are all contained in a ball with volume  $r_n:=\eps \nu_n/2^d$ around $x_1$ and that the WRCM built on $x_1,\dots,x_q$ is fully connected. Here, $\eps>0$ is chosen such that $\vp(\eps)>0$ and such that for all $x \in W_n^{-r_n}$ and all $n \ge 1$, the ball centred at $x$ with volume $r_n$ is entirely contained in $W_n$. Such an $\eps$ exists since the sequence $(\nu_n/n)_n$ is bounded. Since any two points in a ball with volume $r_n$ are connected with probability at least $p:= \vp(\eps)$ and $\rho_q(\xx)\ge K_0(o)^q$ for $\alpha >0$, this yields the lower bound
\begin{align*}
	\int_{\mathbb X^q}  f_n(\mathbf s)  \rho_{q}(\xx) \mu_n^q(\d \mathbf s)  \ge K_0(o)^q|W_n^{-r_n}| p^{\binom q2} r_n^{q-1}.
\end{align*}
With an analogous argument for the term with $k=q-1$ in \eqref{eq:varineq} we find that
\begin{align*}
	\int_{\mathbb X^q} 	\int_{\mathbb X^{q-1}} f_n(\mathbf s) f_n(s_1,\mathbf t) \rho_{2q-1}(\xx,\yy) \mu_n^{q-1}(\d \mathbf t) \mu_n^q(\d \mathbf s)  \ge K_0(o)^{2q-1}|W_n^{-r_n}| p^{2\binom q2} r_n^{2q-2}.
\end{align*}
Combining these estimates gives the assertion.
\enp

\begin{remark} \label{rem:det}
	It is a notoriously challenging task to obtain lower bounds on the variance of geometric functionals of (stationary) determinantal point processes (see also the discussion in Remark (iii) after Theorem 1.14 in \cite{BYY}). This is due to the fact that $\rho_{2q}(\xx,\yy)\le \rho_q(\xx) \rho_q(\yy)$ for $\a\le 0$, which results in positive and negative terms in the expansion of the variance in \eqref{eq:vareq}. In the case $q=1$ (where $S_n$ is the number of points in $\PP \cap W_n$) we have for $\int |K_0(x)|^2 \d x =K_0(o)$ that the variance of $S_n$ grows for $n \to \infty$ as the surface area of $W_n$, i.e.~$\Var(S_n) \sim |W_n|^{\f{d-1}{d}}$ (see \cite{FH99}, (4.15)). For instance, this is the case in the Ginibre process, and similar phenomena also occur for $\int |K_0(x)|^2 \d x >K_0(o)$. However, if  $\int |K_0(x)|^2 \d x < K_0(o)$ we have that $\Var(S_n)\sim |W_n|$. If $q \ge 2$, the situation is more complex, since the asymptotic behavior of the variance of $S_n(G)$ is determined by the asymptotic properties of $K$, $(\nu_n)_{n \ge 1}$ as well as the profile function $\vp$. A concrete scenario in which we can determine the growth rate of $\Var(S_n(G))$ is when $q=2$, $G$ is the connected graph on two vertices, $\nu_n\to 0$ and the profile function $\vp$ has bounded support. In this case we obtain from \eqref{eq:vareq} that
	\begin{align*}
		&\Var(S_n(G))=\int_{W_n^2} \E[\vp_n(x_1-x_2,U_1,U_2)] \r_2(x_1,x_2) \d(x_1,x_2) \\
		&\quad + 2\int_{W_n^3} \E[\vp_n(x_1-x_2,U_1,U_2) \vp_n(x_1-x_3,U_1,U_3)] \r_3(x_1,x_2,x_3) \d(x_1,x_2,x_3)\\
		&\quad + \int_{W_n^4} \E[\vp_n(x_1-x_2,U_1,U_2) \vp_n(x_3-x_4,U_3,U_4)] (\r_4(x_1,x_2,x_3,x_4) - \r_2(x_1,x_2) \r_2(x_3,x_4)) \d(x_1,x_2,x_3,x_4).
	\end{align*}
	A detailed analysis of the three different terms that proceeds along the lines of the proof of \cite[Theorem 3.2]{klein} shows that the first term grows at rate $|W_n| \nu_n$, whereas the second and the third term grow at a rate that is at most $|W_n| \nu_n^2$. Hence, there is some $v>0$ such that
	$$
	\Var(S_n(G))\sim v|W_n| \nu_n\quad \text{as } n \to \infty.
	$$
\end{remark}

\section*{Acknowledgments}
This work was supported by a visit grant by the Danish Data Science Academy (DDSA-V-2023-009).

\bibliographystyle{abbrv}
\bibliography{./main}
\end{document}